\documentclass[12pt]{article}

\newcommand{\ncm}{\newcommand}
\ncm{\aut}{auto\-mor\-phi\-sm} \ncm{\Inn}{\mbox{\rm Inn($A$)}} 
\ncm{\Ap}{\mbox{$\overline{\rm Inn}(A)$}} \ncm{\Ext}{\mbox{\rm 
Ext}} \ncm{\Ex}{\mbox{\rm Ex}} \ncm{\OExt}{\mbox{\rm OrderExt}} 
\ncm{\AI}{\mbox{\rm AInn($A$)}} \ncm{\HI}{\mbox{\rm HInn($A$)}} 
\ncm{\Aut}{\mbox{\rm Aut($A$)}} \ncm{\Mal}{\mbox{$M_{\alpha}$}} 
\ncm{\Aff}{\mbox{${\rm Aff}$}} \ncm{\id}{\mbox{\rm id}} 
\ncm{\Ker}{\mbox{\rm Ker}} \ncm{\BE}{\begin{eqnarray*}} 
\ncm{\EE}{\end{eqnarray*}} \ncm{\lra}{\mbox{$\longrightarrow$}} 
\ncm{\Hom}{\mbox{\rm Hom}} \ncm{\calU}{{\cal U}} \ncm{\el}{\ell} 
\ncm{\ad}{\mbox{\rm ad}} \ncm{\Alg}{\mbox{\rm Alg}} 
\ncm{\Conv}{\mbox{\rm Conv}} \ncm{\D}{{\cal D}}

\ncm{\cstar}{$C^{*}$-algebra} \ncm{\cstars}{$C^{*}$-algebras} 
\ncm{\ra}{\mbox{$\rightarrow$}} \ncm{\la}{\mbox{$\leftarrow$}} 
\ncm{\hra}{\hookrightarrow} \ncm{\da}{\mbox{$\downarrow$}} 
\ncm{\se}{\mbox{$\searrow$}} \ncm{\al}{\mbox{$\alpha $}} 
\ncm{\del}{\mbox{$\delta$}} \ncm{\supp}{\mbox{\rm supp}} 
\ncm{\Ad}{\mbox{\rm Ad}} \ncm{\CAR}{\mbox{$M_{2^{\infty}}$}} 
\ncm{\ep}{\mbox{$\epsilon > 0$}} \ncm{\mod}{\mbox{\rm mod}} 
\ncm{\Sp}{\mbox{\rm Sp}} \ncm{\ol}{\overline} 
\ncm{\Mninf}{\mbox{$M_{n^{\infty}}$}} \ncm{\MR}{M. R\o{}rdam} 
\ncm{\Range}{\mbox{\rm Range}} 
\ncm{\vo}{}
\ncm{\ch}{}
\ncm{\CMP}{Comm. Math. Phys.} \ncm{\add}{} 
\ncm{\tilsig}{\tilde{\sigma}}

\newtheorem{theo}{Theorem}[section]

\newtheorem{lem}[theo]{Lemma}
\newtheorem{prop}[theo]{Proposition}
\newtheorem{remark}[theo]{Remark}
\newtheorem{definition}[theo]{Definition}
\newtheorem{example}[theo]{Example}

\newenvironment{rem}{\begin{remark} \rm}{\end{remark}}
\newenvironment{pf}{{\it Proof.}}{\vspace{3mm}}

\ncm{\R}{\mbox{\bf R}} \ncm{\Z}{\mbox{\bf Z}} \ncm{\T}{\mbox{\bf 
T}} \ncm{\TT}{\T$^{2}$} \ncm{\N}{\mbox{\bf N}} \ncm{\C}{\mbox{\bf 
C}} 

\title{Examples of one-parameter automorphism groups of UHF algebras}
\bigskip

\author{A. Kishimoto\\
  \small Department of Mathematics, Hokkaido University,
   Sapporo 060, Japan}
\date{\small }
\begin{document}
\maketitle

\section{Introduction} 
B. Blackadar \cite{Bl} constructed for the first time an example 
of a symmetry (or an automorphism of period two) of the CAR 
algebra (or the UHF algebra of type $2^{\infty}$) whose fixed 
point algebra is not AF (or approximately finite-dimensional). 
This was soon extended to produce an example of finite-group 
actions on UHF algebras whose fixed point algebras are not AF 
\cite{BEEK1} and then of compact-group actions \cite{EK}. Note 
that these examples can now be obtained as corollaries 
\cite{K1,BEEK2} to the classification results for certain amenable 
\cstars\ started by G.A. Elliott \cite{Ell1} and extended by 
himself and many others (see e.g. \cite{DG}). 

In the same spirit we present yet other examples, this time, of 
one-parameter automorphism groups of UHF algebras, which do not 
seem to follow as a consequence from the above general results. 
Before stating what kind of examples they are we first recall the 
subject from \cite{BR,Br,Sak}. 

Let $A$ be a UHF algebra (or more generally, a simple AF \cstar) 
and let $\alpha$ be a one-parameter automorphism group of $A$. We 
always assume that $t\mapsto \alpha_t(x)$ is continuous for each 
$x\in A$ and denote by $\delta=\delta_{\alpha}$ the 
(infinitesimal) generator of $\alpha$: $$ 
 \delta(x)=\lim_{t\ra0}(\alpha_t(x)-x)/t,\ \ x\in \D(\delta),
 $$ 
where the domain $\D(\delta)$ of $\delta$ is the set of $x\in A$ 
for which the limit exists. Then $\D(\delta)$ is a dense 
*-subalgebra of $A$ and $\delta$ is a *-derivation of $\D(\delta)$ 
into $A$. We equip $\D(\delta)$ with the norm 
$\|\,\cdot\,\|_{\delta}$: 
$$ 
 \|x\|_{\delta}=\|\left(\begin{array}{cc}x &\delta(x)\\ 0 & x
 \end{array}\right)\|_{M_2\otimes A},
 $$
by which $\D(\delta)$ is a Banach *-algebra. We note that the 
\cstar\ $A$ can be recovered as the universal \cstar\ of 
$\D(\delta)$ and that the self-adjoint part of $\D(\delta)$ is 
closed under $C^{\infty}$ functional calculus.  Let us say that a 
Banach *-algebra (or a \cstar) $D$ is AF if $D$ has an increasing 
sequence $(D_n)$ of finite-dimensional *-subalgebras such that the 
union $\cup_n D_n$ is dense in $D$. 

If $\alpha_t$ is given as $\Ad\,e^{ith}$ with some $h=h^*\in A$, 
then $\alpha$ is called inner and the generator $\delta_{\alpha}$ 
is defined on the whole of $A$ and is given by $\ad\,ih$. If 
$\alpha$ is obtained as the limit (pointwise on $A$ and uniformly 
on compact subsets of $\R$) of inner one-parameter automorphism 
groups, then $\alpha$ is called approximately inner. 

For a one-parameter automorphism group $\alpha$ of the UHF algebra 
$A$ with generator $\delta_{\alpha}$ we quote the following two 
results of S. Sakai \cite{Sak0,Sak}: 

\begin{theo}\label{A} 
The Banach *-algebra $\D(\delta_{\alpha})$ contains an AF Banach 
*-subalgebra $B$ such that $B$ is dense in $A$ under the embedding 
$\D(\delta_{\alpha})\subset A$. \end{theo} 

\begin{theo}\label{B}
If $\D(\delta_{\alpha})$ is AF, then $\alpha$ is approximately 
inner. \end{theo} 
 
The condition of Theorem \ref{B} was satisfied for all the known 
examples so far and the core problem (\cite{Sak}, 4.5.10), as a 
possible solution to the Powers-Sakai conjecture (\cite{Sak}, 
4.5.9), asks whether this is true for all one-parameter 
automorphism groups of UHF algebras. We will give an example which 
shows this is not the case; we construct an approximately inner 
one-parameter automorphism group $\alpha$ such that 
$\D(\delta_{\alpha})$ is not AF (see \ref{C}). The property we use 
to conclude this  is real rank \cite{BP}.  

So far we know of no examples of one-parameter automorphism groups 
$\alpha$ of UHF algebras such that $\D(\delta_{\alpha})$ contains 
no maximal abelian $C^*$-subalgebra (masa) of $A$. But we will 
present another example which shows that $\D(\delta_{\alpha})$ 
need not contain a canonical AF masa of $A$ even though 
$\D(\delta_{\alpha})$ is AF (see \ref{3D}). (We call an abelian 
$C^*$-subalgebra $C$ of $A$ a canonical AF masa if there exists an 
increasing sequence $(A_n)$ of finite-dimensional 
$C^*$-subalgebras of $A$ such that $A=\ol{\cup_nA_n}$ and $C\cap 
A_n\cap A_{n-1}'$ is a masa of $A_n\cap A_{n-1}'$ for all $n$ with 
$A_0=0$.)  As will be shown in \ref{3A}, this is equivalent to the 
property that  any inner perturbation of $\alpha$ is not AF 
locally representable. (We call $\alpha$ {\em AF locally 
representable} if there exists an increasing sequence $(A_n)$ of 
finite-dimensional $C^*$-subalgebras of $A$ with dense union such 
that $\alpha$ leaves each $A_n$ invariant and thus $\alpha|A_n$ is 
inner. In this case there is a canonical AF masa $C$ associated 
with $(A_n)$ such that $\delta_{\alpha}|C=0$, and the union 
$\cup_nA_n$ is a core for the generator $\delta_{\alpha}$ and thus 
$\D(\delta_{\alpha})$ is AF.) In our example, if $\alpha$ is 
periodic, we use the property that the fixed point algebra 
$A^{\alpha}$ is not AF to conclude that $D(\delta_{\alpha})$ 
contains no canonical AF masa. If $\alpha$ is not periodic, we 
instead use the property that for some unitary $u$ in 
$\D(\delta_{\alpha})$ with $\delta_{\alpha}(u)\approx0$ there is 
no continuous path $(u_t)$ of unitaries between $u$ and $1$ such 
that $\delta_{\alpha}(u_t)\approx0$. 

To sum up let us state three properties for $\alpha$: 
 
(1) $\D(\delta_{\alpha})$ contains a canonical AF masa.

(2) $\D(\delta_{\alpha})$ is AF.

(3) $\alpha$ is approximately inner. 

\noindent Then $(1)\Rightarrow(2)\Rightarrow(3)$ but 
$(1)\not\Leftarrow(2)\not\Leftarrow(3)$. 

\section{A generator whose domain is not AF}
 \setcounter{theo}{0} 
\begin{theo}\label{C}
Let $A$ be a non type I simple AF algebra. Then there exists an 
approximately inner one-parameter automorphism group $\alpha$ of 
$A$ such that the domain $\D(\delta_{\alpha})$ is not AF. 
\end{theo} 
 \begin{pf}
Let $(A_n)$ be an increasing sequence of finite-dimensional 
$C^*$-subalgebras of $A$ such that $A=\ol{\cup_nA_n}$ and let 
$A_n=\oplus_{j=1}^{k_n}A_{nj}$ be the direct sum decomposition of 
$A_n$ into full matrix algebras $A_{nj}$. Since 
$K_0(A_n)\cong\Z^{k_n}$, we obtain a sequence of $K_0$ groups: $$
 \Z^{k_1}\stackrel{\chi_1}{\ra}\Z^{k_2}\stackrel{\chi_2}{\ra}\cdots,
 $$
where $\chi_n$ is the positive map of $K_0(A_n)=\Z^{k_n}$ into 
$A_{n+1}=\Z^{k_{n+1}}$ induced by the embedding $A_n\subset 
A_{n+1}$. Since $K_0(A)$ is a simple dimension group other than 
$\Z$, we may assume that all $\chi_n(i,j)\geq 3$. 

By using $(A_n)$ we will express $A$ as an inductive limit of 
\cstars\ $A_n\otimes C[0,1]$. First we define a homomorphism 
$\varphi_{n,ij}$ of $A_{nj}\otimes C[0,1]$ into $A_{nj}\otimes 
M_{\chi_n(i,j)}\otimes C[0,1]$, with $M_k$ the full $k$ by $k$ 
matrix algebra, as follows:
 $$
 \varphi_{n,ij}(x)(t)=x(t)\oplus 
 \oplus_{\ell=0}^{\chi_n(i,j)-2}x(\frac{t+\ell}{\chi_n(i,j)-1}),
 $$
in particular $\varphi_{n,ij}(x)$ is of diagonal form in the 
matrix algebra over $A_{nj}\otimes C[0,1]$. (In the above 
definition of $\varphi_{n,ij}$, the variable $t$ inside the 
summands after the first could be removed if we assumed that 
$\min_{ij}\chi_n(i,j)\ra\infty$.) Then on the embedding  of
 $$
 \oplus_{j=1}^{k_n}A_{nj}\otimes M_{\chi_n(i,j)}\otimes C[0,1]
 $$
in the natural way into $A_{n+1,i}\otimes C[0,1]$, 
$(\varphi_{n,ij})$ defines an injective homomorphism 
$\varphi_n:\,A_n\otimes C[0,1]\ra A_{n+1}\otimes C[0,1]$. Then it 
follows that the inductive limit of the sequence $(A_n\otimes 
C[0,1],\varphi_n)$ is isomorphic (as we shall explain) to the 
given \cstar\ $A$; we have thus expressed $A$ as $\ol{\cup_nB_n}$ 
where $B_n=A_n\otimes C[0,1]\subset B_{n+1}$. 

This isomorphism follows from Elliott's result \cite{Ell1} by 
checking that the inductive limit \cstar\ is simple and of real 
rank zero \cite{BBEK} and has the right K-theoretic data, these 
properties being consequences of the condition $\chi_n(i,j)\geq 3$ 
and the special form of $\varphi_{n,ij}$. (As a matter of fact, 
the proof that the inductive limit is AF is about the same as the 
proof that it has real rank zero, since both of these facts follow 
by showing that the canonical self-adjoint element $x_n\in 
1\otimes C[0,1]\subset B_n$, see below, can be approximated by 
self-adjoint elements with finite spectra in $A\cap (A_n\otimes 
1)'$.) 

We will define a one-parameter automorphism group $\alpha$ of $A$ 
with the property $\alpha_t(B_n)=B_n$. First we define a sequence 
$(H_n)$ with self-adjoint $H_n\in A_n\otimes 1\subset B_n$ 
inductively. Let $H_1\in A_1\otimes 1\subset B_1$ and let 
$H_{n}=H_{n-1}+\sum_i\sum_j h_{n,ij}$, where $h_{n,ij}$ is a 
self-adjoint diagonal matrix of $M_{\chi_{n-1}(i,j)}$ which is
identified with a $C^*$-subalgebra of $B_n$ by 
 $$
 M_{\chi_{n-1}(i,j)}\equiv 1\otimes M_{\chi_{n-1}(i,j)}\otimes 1
 \subset A_{n-1,j}\otimes M_{\chi_{n-1}(i,j)}\otimes 1\subset B_n.
 $$
We define $\alpha_t|B_n$ by $\Ad\,e^{itH_n}|B_n$. Since 
$\alpha_t|B_n=\Ad\,e^{itH_{n+1}}|B_n$ from the definition of 
$H_{n+1}$, $(\alpha_t|B_n)$ defines a one-parameter automorphism 
group $\alpha$ of $A$. Let $x$ be the identity function on the 
interval $[0,1]$ and let $x_n=1\otimes x\in 1\otimes C[0,1]\subset 
B_n$. Then we have that $\alpha_t(x_n)=x_n,\ t\in \R$, or 
$\delta(x_n)=0$ for the generator $\delta$ of $\alpha$.

Note that $\D(\delta)$ contains $\cup_nB_n$. Since $(1\pm 
\delta)\cup_nB_n=\cup_nB_n$ and $\|(1\pm\delta)x\|\geq \|x\|$ for 
$x\in\D(\delta)$, it follows that $\cup_nB_n$ is a core for 
$\delta$, i.e., $\cup_nB_n$ is dense in the Banach *-algebra 
$\D(\delta)$. See \cite{BR,Sak} for details.

If $\D(\delta)$ is AF, then for any $h\in \D(\delta)_{sa}=\{y\in 
D(\delta)\,; y=y^*\}$ there  exists a sequence $(h_n)$ in 
$\D(\delta)_{sa}$ such that $\Sp(h_n)$ is finite and 
$\|h-h_n\|_{\delta}\ra0$. (We can further impose, without 
difficulty, the condition on $(h_n)$ that $\Sp(h_n)$ is a subset 
of the smallest closed interval containing $\Sp(h)$.) Here we note 
that the spectrum of $h_n$, $\Sp(h_n)$, may be computed in 
$\D(\delta)$ or in $A$ since they are the same. (If $y\in 
\D(\delta)$ is invertible in $A$, then it follows that $y^{-1}\in 
\D(\delta)$ and $\delta(y^{-1})=-y^{-1}\delta(y)y^{-1}$.) In this 
case we may say that the Banach *-algebra $\D(\delta)$ has real 
rank zero as in the case of \cstars\ \cite{BP}. What we will do is  
show that $\D(\delta)$ does not have real rank zero for certain 
$\alpha$ (and hence has real rank one, by defining real rank for 
$\D(\delta)$ as in \cite{BP}). 

We fix $H_1$ and all $h_{n,ij}$ except for $h_{n,11}$. We will 
inductively define $h_{n,11}$ to be of the form 
 $$
 a_n\oplus0\oplus \cdots \oplus0\in 1\otimes M_{\chi_{n-1}(1,1)}\otimes 1
\subset A_{n1}\otimes C[0,1]
 $$
with $a_n>0$, to make sure that no $x_n$ can be approximated by 
self-adjoint elements with finite spectra in $\D(\delta)$. 
                                          
Let $P_n$ be the identity of $B_{n1}=A_{n1}\otimes C[0,1]$ and let 
$Q_n$ be the projection
 $$
 1\oplus0\oplus\cdots\oplus0\in 1\otimes M_{\chi_{n-1}(1,1)}\otimes1\subset B_{n1}.
 $$
Let $(\epsilon_m)$ be a strictly decreasing sequence of positive 
numbers such that $\epsilon_1\leq 3/5$. We shall construct a 
sequence $(a_n)$  such that if $h=h^*\in B_{m+1,1}$ satisfies that 
$0\leq h\leq 1,\ \mu(\Sp(h))<\epsilon_{m+1}$, and 
$\|h-x_nP_{m+1}\|<1/5$ for some $n\leq m$, then $\|\delta(h)\|>1$, 
where $\mu$ denotes  Lebesgue measure on \R. (Here we have imposed 
the condition $0\leq h\leq1$, which does not cause any loss of 
generality.) This in particular shows that if $h=h^*\in A$ with 
$0\leq h\leq1$ belongs to $\cup_n B_n$ and satisfies that $\Sp(h)$ 
is finite and $\|x_n-h\|<1/5$ for some $n$, then 
$\|\delta(h)\|>1$. We will discuss later how to remove the 
condition $h\in \cup_nB_n$ in this statement.  

Let $m=1$ and choose $a_1$ arbitrarily. If $h\in B_{2,1}$ with 
$0\leq h\leq 1,\ \mu(\Sp(h))<\epsilon_2$, and 
$\|\delta(h)\|\leq1$, then 
 \BE
 \|\delta(h)\|&\geq& \|Q_2[iH_1+\sum_i\sum_jih_{2,ij},h](1-Q_2)\|\\
 &\geq& a_2\|Q_2h(1-Q_2)\|-2\|H_1\|-\sum_{j=2}^{k_1}\|h_{2,1j}\|
 \EE
since $Q_2h_{2,11}=a_2Q_2$ and $Q_2h_{2,ij}=0$ for 
$(i,j)\neq(1,1)$. If $a_2$ is sufficiently large, then 
 $$
 \|Q_2h(1-Q_2)\|<\frac{\epsilon_1-\epsilon_2}{4\xi_2},
 $$
where $\xi_n$ is defined by $A_{n,1}\cong M_{\xi_n}$. If 
$\tilde{h}=Q_2hQ_2+(1-Q_2)h(1-Q_2)$, then we have  that 
 $$
 \|h-\tilde{h}\|< 
 \frac{\epsilon_1-\epsilon_2}{2\xi_2}.
 $$
Since $h\in A_{2,1}\otimes C[0,1]$, there are at most $\xi_2$ 
eigenvalues of $h(t)$ for each $t\in [0,1]$. Thus $\Sp(h)$ 
consists of at most $\xi_2$ closed intervals. Since 
 $$
 \Sp(\tilde{h})\subset\Sp(h)+[-\|h-\tilde{h}\|,\|h-\tilde{h}\|],
 $$
we obtain that
 $$
 \mu(\Sp(\tilde{h}))<\epsilon_2 +2\xi_2\|h-\tilde{h}\|<\epsilon_1\leq 3/5,
 $$
which, in particular, implies that $\mu(\Sp(Q_2hQ_2))<3/5$. On the 
other hand, if $\|x_1P_2-h\|<1/5$, then $\|x_1Q_2-Q_2hQ_2\|<1/5$, 
which implies that
 \BE
 \|Q_2h(0)Q_2\|&<&1/5\\
 \|Q_2-Q_2h(1)Q_2\|&<&1/5.
 \EE
Since $t\mapsto \Sp(Q_2h(t)Q_2)$ is continuous in a certain 
well-known sense, it follows that $\Sp(Q_2hQ_2)$ strictly contains 
$[1/5,4/5]$, and so $\mu(\Sp(Q_2hQ_2))>3/5$. This contradiction 
completes the proof for the case $m=1$. 

Suppose that $a_1,\ldots, a_m$ are chosen in such a way that the 
conditions for $h\in B_{n+1,1}$ with $n<m$ are satisfied. Let 
$h\in B_{m+1,1}$ be such that $0\leq h\leq 1$, 
$\mu(\Sp(h))<\epsilon_{m+1}$, and $\|\delta(h)\|\leq1$. Then as 
before we have that 
 \BE
 \|\delta(h)\|&\geq& \|Q_{m+1}\delta(h)(1-Q_{m+1})\|\\
 &\geq&a_{m+1}\|Q_{m+1}h(1-Q_{m+1})\|-2\|H_m\|-\sum_{j=2}^{k_{m+1}}\|h_{m+1,1j}\|.
 \EE
We choose $a_{m+1}$ in such a way that $\|\delta(h)\|\leq1$ 
implies that
 $$
 \|Q_{m+1}h(1-Q_{m+1})\|<\frac{\epsilon_m-\epsilon_{m+1}}{4\xi_{m+1}}.
 $$
Then for $\tilde{h}=Q_{m+1}hQ_{m+1}+(1-Q_{m+1})h(1-Q_{m+1})$ we 
have that 
 $$
 \|h-\tilde{h}\|< 
 \frac{\epsilon_m-\epsilon_{m+1}}{2\xi_{m+1}}.
 $$
Since $\Sp(h)$ consists of at most $\xi_{m+1}$ connected 
components and $\mu(\Sp(h))<\epsilon_{m+1}$, it follows that 
$\mu(\Sp(\tilde{h}))<\epsilon_m$. If $\|x_nP_{m+1}-h\|<1/5$ for 
some $n<m+1$, then 
 $$
 \|Q_{m+1}x_n-Q_{m+1}hQ_{m+1}\|<1/5.
 $$
If $n=m$, then by using that $Q_{m+1}x_m(0)=0$ and 
$Q_{m+1}x_m(1)=Q_{m+1}$ we can reach a contradiction as before. If 
$n<m$, then since $\varphi_m$ restricts to an isomorphism  of 
$B_{m1}=A_{m1}\otimes C[0,1]$ onto $Q_{m+1}B_{m+1,1}Q_{m+1}$ 
mapping $x_nP_m$ into $x_nQ_{m+1}$, the pre-image $k \in B_{m1}$ 
of $Q_{m+1}hQ_{m+1}$ satisfies that $0\leq k\leq 1,\ 
\mu(\Sp(k))<\epsilon_m$, $\|\delta(k)\|\leq1$, and 
$\|x_nP_m-k\|<1/5$. Thus this would be a contradiction by the 
induction hypothesis.
                                                  
Thus we have shown that if $h=h^*\in\cup_nB_n$ has finite spectrum 
in $[0,1]$ and $\|h-x_n\|<1/5$ for some $n$ then  
$\|\delta(h)\|>1$. 

Let $h=h^*\in \D(\delta)$ be such that $\Sp(h)$ is a finite subset 
of $[0,1]$. Since $\cup_nB_n$ is a core for $\delta$ (or 
$\cup_nB_n$ is dense in $\D(\delta)$), we obtain a sequence 
$(h_n)$ in $\cup_nB_n$ such that $\|h-h_n\|\ra0$ and 
$\|\delta(h)-\delta(h_n)\|\ra0$. Obviously we may suppose that 
$h_n=h_n^*$. If $\Sp(h)=\{\lambda_1,\lambda_2,\ldots,\lambda_k\}$ 
and $\epsilon>0$ is smaller than any $|\lambda_i-\lambda_j|/2$ 
with any $i\neq j$, $\Sp(h_n)$ is covered by the disjoint union of 
the $\epsilon$-neighborhoods of 
$\lambda_1,\lambda_2,\ldots,\lambda_k$ for all sufficiently large 
$n$. If we define 
 $$
 p_{ni}=\frac{1}{2\pi i}\int_{|z-\lambda_i|=\epsilon}(z-h_n)^{-1}dz,
 $$
which is a projection in $\cup_nB_n$, we have that $ 
\|h_n-\sum_i\lambda_ip_{ni}\|\approx0$, depending on 
$\|h_n-h\|\approx0$. If we define a projection $p_i$ just as 
$p_{ni}$ by using $h$ instead of $h_n$, we obtain that 
$p_i\in\D(\delta)$, $h=\sum_i\lambda_ip_i$, $\|p_{ni}-p_i\|\ra0$, 
and $\|\delta(p_{ni})-\delta(p_i)\|\ra0$. (Here the latter 
convergence follows by using 
$\delta((z-h_n)^{-1})=(z-h_n)^{-1}\delta(h_n)(z-h_n)^{-1}$.) Hence 
it follows that $\tilde{h}_n=\sum_i\lambda_ip_{ni}\in \cup_nB_n$ 
satisfies that $\Sp(\tilde{h}_n)=\Sp(h)$, $\|\tilde{h}_n-h\|\ra0$, 
and $\|\delta(\tilde{h}_n)-\delta(h)\|\ra0$. In this way if 
furthermore $\|h-x_n\|<1/5$ for some $n$, we can conclude that 
$\|\delta(h)\|>1$. (Here to get the strict inequality instead of 
$\|\delta(h)\|\geq1$, we may apply this argument to $\lambda h$ 
with $0<\lambda<1$ and $\|\lambda h-x_n\|<1/5$ instead of the 
given $h$.) 
    
\end{pf}
   
In the situation of the above proof we define a masa (maximal 
abelian $C^*$-subalgebra) $C_n$ of $A_n$ with $H_n\in C_n$ as 
follows. Let $C_1$ be a  masa of $A_1$ containing $H_1$. We 
inductively define a masa $C_n$ of $A_n$ as the $C^*$-subalgebra 
generated by $C_{n-1}$ and a masa $A_n\cap A_{n-1}'$ containing 
$h_{n-1,ij}$ for all $i,j$. Then $C=\ol{\cup_nC_n}$ is a maximal 
abelian AF $C^*$-subalgebra of $\tilde{A}=\ol{\cup_nA_n}$, which 
is regarded as a $C^*$-subalgebra of $A=\ol{\cup_nB_n}$ and is 
isomorphic to $A$ itself. Let $D_n$ be the $C^*$-subalgebra of 
$B_n=A_n\otimes C[0,1]$ generated by $C_n$ and $1\otimes C[0,1]$. 
Then $(D_n)$ forms an increasing sequence and generates a masa $D$ 
of $A$. Since our generator $\delta$ vanishes on $D$, $\D(\delta)$ 
contains the $C^*$-algebra $D$. Furthermore $D$ is Cartan in the 
sense that the unitary subgroup $\{u\in A\,;\, uDu^*=D\}$ 
generates $A$ while $C$ is Cartan in $\tilde{A}$.  Since $\alpha$ 
fixes $\tilde{A}$, we may consider the one-parameter automorphism 
group $\tilde{\alpha}=\alpha|\tilde{A}$, whose generator 
$\tilde{\delta}$ vanishes on $C$. Since $\alpha$ (resp. 
$\tilde{\alpha}$) fixes the generating sequence $(B_n)$ of $A$ 
(resp. $(A_n)$ of $\tilde{A}$) and is inner on each $B_n$ (resp. 
$A_n$), both $\alpha$ and $\tilde{\alpha}$ can be called {\em 
locally representable} (or locally inner). (To be more specific, 
we call a one-parameter automorphism group $\alpha$ of a 
$C^*$-algebra $A$ locally representable if there is an increasing 
sequence $(A_n)$ of $C^*$-subalgebras of $A$ such that $\alpha$ 
leaves $A_n$ invariant and $\alpha|A_n$ is inner for each $n$ and 
the union $\cup_nA_n$ is dense in $A$.) We will call 
$\tilde{\alpha}$ especially {\em AF locally representable} (or AF 
representable) since the $A_n$'s are finite-dimensional. (What we 
meant by locally representable in \cite{K5} is in this latter 
sense.) Thus $\alpha$ does not look very queer as $\tilde{\alpha}$ 
does not, which may make the following problem interesting: 

\medskip
\noindent {\it Problem} For a unital simple AF \cstar\ $A$ is 
there a one-parameter automorphism group $\alpha$ of $A$ such that 
$\alpha$ is not an inner perturbation of a locally representable 
one? 

\medskip
As a matter of fact this is probably what is meant by \cite{Sak}, 
4.5.8. A locally representable one-parameter automorphism group 
$\alpha$  might be characterized, up to inner perturbation, by the 
property that $\D(\delta_{\alpha})$ contains a masa of $A$. (In 
any case a similar problem is to find $\alpha$ such that  
$\D(\delta_{\alpha})$ contains no masa.) Though this looks 
optimistic, we will consider a special case in the next section. 

We shall conclude this section with a remark on {\em commutative 
normal *-derivations} introduced by S. Sakai \cite{Sak}.

A commutative normal *-derivation $\delta$ in an AF \cstar\ $A$ is 
defined as follows: $\delta$ is defined on $\D(\delta)=\cup_nA_n$ 
for some increasing sequence $(A_n)$ of finite-dimensional 
$C^*$-subalgebras of $A$ with dense union and has a mutually 
commuting family $\{h_n\}$ of self-adjoint elements in $A$ such 
that $\delta|A_n=\ad\,ih_n|A_n$. (This is adapted from 4.1.5 and 
4.5.7 in \cite{Sak} to the case of AF \cstars.) Then it follows 
that $\delta$ extends to a generator, which generates an 
approximately inner one-parameter automorphism group (see 
\cite{Sak}, 4.1.11, and \cite{BK} for a similar result). 

\begin{rem}
If $A$ is a non type I simple AF \cstar, there is a commutative 
normal *-derivation whose closure is not a generator. 
\end{rem}

Such an example is given in the proof of Theorem \ref{C}. The 
proof that $\ol{\cup_nB_n}$ with $B_n=A_n\otimes C[0,1]$ is AF (or 
even just the fact that this \cstar\ is AF)  shows that any finite 
subset in $\cup_nB_n$ can be approximately contained in a 
finite-dimensional $C^*$-subalgebra of $B_m$ for some $m$. Hence 
we can construct an increasing sequence $(D_n)$ of 
finite-dimensional $C^*$-subalgebras in $\cup_nB_n$ such that 
$\ol{\cup_nD_n}=\ol{\cup_nB_n}$. Thus there is a subsequence 
$(k_n)$ such that $D_n\subset B_{k_n}$, which shows that 
$\delta|D_n=\ad\,iH_{k_n}|D_n$. This way 
$\delta_0=\delta|\cup_nD_n$ is a commutative *-derivation. If the 
closure $\ol{\delta_0}$ were a generator, then it must be $\delta$ 
and $\D(\ol{\delta_0})$ would be AF, which is a contradiction.

\section{AF locally representable actions} \setcounter{theo}{0}

When $A$ is a unital simple AF \cstar\ and $C$ is a maximal 
abelian AF $C^*$-subalgebra of $A$, we call $C$ a canonical AF 
masa if there is an increasing sequence $(A_n)$ of 
finite-dimensional $C^*$-subalgebras of $A$ such that $\cup_nA_n$ 
is dense in $A$ and $C\cap A_n\cap A_{n-1}'$ is a masa of $A_n\cap 
A_{n-1}'$ with $A_0=0$. Hence $C$ is generated by $C_n=C\cap 
A_n\cap A_{n-1}',\ n=1,2,\ldots$; there is a natural homomorphism 
from the infinite tensor product $\otimes_{n=1}^{\infty}C_n$ onto 
$C$. (See \cite{SV} for this kind of masa.) Then we note the 
following: 

\begin{prop}\label{3A}
Let $\alpha$ be a one-parameter automorphism group of a unital 
simple AF \cstar\ $A$. Then the following conditions are 
equivalent: 
 \begin{enumerate}
 \item $\D(\delta_{\alpha})$ contains a canonical AF masa of $A$.
 \item There is an $h=h^*\in A$ and an increasing sequence $(A_n)$ of 
finite-dimensional $C^*$-subalgebras of $A$ such that $\cup_nA_n$ 
is dense in $\D(\delta_{\alpha})=\D(\delta_{\alpha}+\ad\,ih)$ and 
$\delta_{\alpha}+\ad\,ih$ leaves $A_n$ invariant, i.e., an inner 
perturbation of $\delta_{\alpha}$ generates an AF locally 
representable one-parameter automorphism group of $A$. 
\end{enumerate}
\end{prop} 

We only have to show (1) implies (2).  When $C$ denotes the 
canonical masa contained in $\D(\delta)$ with 
$\delta=\delta_{\alpha}$ and $(A_n)$ denotes the associated 
increasing sequence, the proof will go as follows. We first find a 
self-adjoint $h\in A$ such that $\delta|C=-\ad\,ih|C$. Then 
$\delta+\ad\,ih$ vanishes on $C$ and we may take $\delta+\ad\,ih$ 
for $\delta$ and assume that $\delta|C=0$. We then modify $(A_n)$ 
by employing a method in \cite{Sak0} in such a way that 
$\cup_nA_n\subset \D(\delta)$. Next we find a self-adjoint $h_n\in 
C$ such that $\delta|A_n=\ad\,ih_n|A_n$. Using the fact that $h_n$ 
can be chosen from $C$ we find a self-adjoint $h\in C$ such that 
$\delta+\ad\,ih$ leaves $\cup_nA_n$ invariant \cite{Sak}. Thus by 
taking $\delta+\ad\,ih$ for $\delta$ and passing to a subsequence 
of $(A_n)$, we may assume that $C\cup (\cup_nA_n) \subset 
\D(\delta),\ \delta|C=0$, and $\delta(A_n)\subset A_{n+1}$. Then 
we find a self-adjoint $h_n\in C\cap A_{n+1}$ such that 
$\delta|A_n=\ad\,ih_n|A_n$. If we denote by $B_n$ the 
$C^*$-subalgebra generated by $A_n$ and $C\cap A_{n+1}$, then 
$(B_n)$ is an increasing sequence of finite-dimensional 
$C^*$-subalgebras of $A$ with $\ol{\cup_nB_n}=A$ and 
$\delta(B_n)\subset B_n$ for all $n$. Thus this will complete the 
proof. 

In the above argument most of the steps are either straightforward 
or given in \cite{Sak}. An exception may be the assertion made in 
the very beginning, which we shall show below. 

\begin{lem}\label{3B}
There exists an $h=h^*\in A$ such that $\delta(x)=\ad\,ih(x),\ 
x\in C$. 
\end{lem}
\begin{pf}
Since $\D(\delta)\supset C$ and $C$ is a \cstar, $\delta|C$ is 
bounded \cite{Rin}. We first show that 
 $$
 \|\delta|C\cap A_n'\|\ra0.
 $$
Suppose, on the contrary, that there is an $\epsilon>0$ such that
$\|\delta|C\cap A_n'\|>\epsilon$ for all $n$. Since the closure of 
the convex hull of the projections ${\cal P}(C\cap A_n')$ in 
$C\cap A_n'$ equals $\{h\in C\cap A_n'\ ;\ 0\leq h\leq 1\}$, we 
may assume that there is an $\epsilon>0$ such that $\|\delta|{\cal 
P}(C\cap A_n')\|>\epsilon$ for all $n$. Thus we can find a 
sequence $(e_n)$ of projections with $e_n\in C\cap A_n'$ such that 
$\|\delta(e_n)\|>\epsilon$. Since 
$\delta(e_n)=e_n\delta(e_n)(1-e_n)+(1-e_n)\delta(e_n)e_n$ and 
$\|e_n\delta(e_n)(1-e_n)\|>\epsilon$, there exists a state 
$\varphi_n$ of $A$ for any $\gamma\in (0,1)$ such that 
 $$
 \varphi_n(e_n)=1-\gamma,\ \ 
 \varphi_n\delta(e_n)>2\epsilon\sqrt{\gamma(1-\gamma)}.
 $$ 
By using this fact and an approximation argument, we can see that 
for a subsequence $(k_1,k_2,\ldots,k_n)$ the norm of
 \BE
 \delta(e_{k_1}e_{k_2}\cdots e_{k_n})&=&\delta(e_{k_1})e_{k_2}\cdots 
 e_{k_n}\\
 &+&e_{k_1}\delta(e_{k_2})e_{k_3}\cdots e_{k_n}\\
 &+&\cdots\\
 &+&e_{k_1}\cdots e_{k_{n-1}}\delta(e_{k_n})
 \EE
exceeds 
$2n\epsilon\sqrt{1/n(1-1/n)}(1-1/n)^{n-1}\approx2\epsilon\sqrt{n}e^{-1}$, 
since the products almost become {\em tensor products} in the 
above equality. Here we use the fact that $A$ is simple. This is a 
contradiction for a large $n$.   
   
Thus we have shown that $\|\delta|C\cap A_n'\|\ra0$. By passing to 
a subsequence of $(A_n)$ we may suppose that 
 $$
 \sum_n \|\delta|C\cap A_n'\|<\infty.
 $$
Denoting by $G_n$ the unitary group of $C\cap A_n\cap A_{n-1}'$ 
with $A_0=0$, we consider the following integral with respect to 
normalized Haar measures: 
 \BE
 ih_n&=&\int_{G_1\times\cdots\times G_n}\delta(g_1^*g_2^*\cdots g_n^*)g_n\cdots g_1 
 \\
     &=&\int_{G_1}\delta(g_1^*)g_1+\int_{G_1\times 
     G_2}g_1^*\delta(g_2^*)g_2g_1+\cdots \\
     &&+\int_{G_1\times\cdots\times G_n}g_1^*\cdots 
     g_{n-1}^*\delta(g_n^*)g_ng_{n-1}\cdots g_1,
 \EE
which we can see converges as $n\ra\infty$. Since $[ih_n, 
g]=\delta(g)$ for any unitary $g$ in $G_1\times\cdots\times G_n$, 
it follows that $\delta|C\cap A_{n}=\ad\,ih_n|C\cap A_{n}$, which 
completes the proof. 
\end{pf}

In the above proof of (1)$\Rightarrow$(2) we did not really use 
the fact that $\delta_{\alpha}$ is a generator; so we have: 

\begin{rem}\label{3C}
If $\delta$ is a closed *-derivation in an AF \cstar\ $A$ and 
$\D(\delta)$ contains a canonical AF masa of $A$, then $\delta$ is 
a generator and $\D(\delta)$ is AF. \end{rem}

\begin{prop}\label{3D}
Let $A$ be a non type I simple AF \cstar. Then there exists a 
one-parameter automorphism group $\alpha$ of $A$ such that 
$\D(\delta_{\alpha})$ is AF but does not contain a canonical AF 
masa of $A$. \end{prop} 

If $A$ is a UHF algebra of type $p^{\infty}$ for some $p>1$, 
examples for $\T=\R/\Z$ from \cite{EK} will be the desired ones. 

To deal with a simple AF \cstar\ $A$ let us proceed as follows. 
First express $A$ as $\ol{\cup_nA_n}$ with $A_n$ 
finite-dimensional, as in the proof of \ref{C}. With the notation 
there, we assume this time that all the multiplicities 
$\chi_n(i,j)\geq 4$. We define a homomorphism $\varphi_{n,ij}$ of 
$A_{nj}\otimes C(\T)$ into $A_{nj}\otimes M_{\chi_n(i,j)}\otimes 
C(\T)$ by 
 $$
 \varphi_{n,ij}(x)(z)=\left(\begin{array}{cc}0&z\\1&0\end{array}\right)\oplus
  \left(\begin{array}{cc}0&\ol{z}\\1&0\end{array}\right)\oplus 
  \oplus_{\ell=0}^{\chi_n(i,j)-5}x(1),
 $$
and define accordingly $\varphi_n:A_n\otimes C(\T)\ra 
A_{n+1}\otimes C(\T)$. Then it follows \cite{Ell1} that the 
inductive limit \cstar\ of $(A_n\otimes C(\T),\varphi_n)$ is 
isomorphic to the original $A$; we have thus expressed $A$ as 
$\ol{\cup_nB_n}$ where $B_n=A_n\otimes C(\T)\subset B_{n+1}$. 

We define a sequence $(H_n)$ with self-adjoint $H_n\in A_n\otimes 
1\subset B_n$ by $H_1=0$ and $H_n=H_{n-1}+\sum_i\sum_j h_{n,ij}$, 
where $h_{n,ij}\in 1\otimes M_{\chi_{n-1}(i,j)}\otimes 1\subset 
B_n$ is given by 
 $$
 h_{n,ij}=\left(\begin{array}{cc}1&0\\0&1\end{array}\right)\oplus  
          \left(\begin{array}{cc}-1&0\\0&-1\end{array}\right)\oplus0\oplus\cdots\oplus0.
 $$
Then we define a one-parameter automorphism group $\alpha$ of $A$ 
by $\alpha_t|B_n=\Ad\,e^{itH_n}|B_n$. Note that $\Sp(H_n)\subset 
\Z$ and $\alpha_{2\pi}=\id$. Then we can easily conclude that the 
fixed point algebra $A^{\alpha}$ is not AF; $K_1(A^{\alpha})$ is 
not trivial. What we need is this property to conclude that 
$\D(\delta_{\alpha})$ does not contain a canonical AF masa of $A$. 
Before proving this as a lemma below we shall show that 
$\D(\delta_{\alpha})$ is AF. 

Let $z$ be the canonical unitary in $C(\T)$ and let $z_n=1\otimes 
z\in A_n\otimes C(\T)=B_n$.  We have to use an estimate in the 
approximation of $z_n$ by a unitary with finite spectrum in $B_m$ 
for $m>n$. Let $u$ be the image of $z_n$ in $B_m$. Then the part 
of $u(z)$ which is not constant in $z=e^{2\pi t},\ t\in [0,1)$, 
has an equal number of  eigenvalues $\exp(\pm i2\pi 2^{n-m}(t+k))$ 
with $k=0,1,\ldots,2^{m-n}-1$. By using this we approximate $u$ by 
a unitary $v\in B_m\cap (A_n\otimes 1)'$ with finite spectrum with 
the order $\|u-v\|\approx 2^{n-m}$ (see \cite{BBEK}). But the norm 
of $\delta_{\alpha}|B_m\cap (A_n\otimes 1)'$ can be estimated as 
$m-n$, which yields 
$\|\delta_{\alpha}(u)-\delta_{\alpha}(v)\|\leq(m-n)\|u-v\|$. Thus 
we can conclude that we can make the approximation in 
$\|\,\cdot\,\|_{\delta_{\alpha}}$, which shows that 
$\D(\delta_{\alpha})$ is AF.  

\begin{lem}\label{3E}
If $\alpha$ is a periodic one-parameter automorphism group of a 
simple AF \cstar\ $A$ and $A^{\alpha}$ is not AF, then 
$\D(\delta_{\alpha})$ does not contain a canonical AF masa of $A$.
\end{lem}
\begin{pf}
Let $\delta=\delta_{\alpha}$ and suppose that $\D(\delta)$ 
contains a canonical AF masa $C$. Then by \ref{3A} we have a 
self-adjoint $h\in A$ and an increasing sequence $(A_n)$ of 
finite-dimensional $C^*$-subalgebras of $A$ with dense union such 
that $\delta+\ad\,ih$ leaves $A_n$ invariant and $C\cap A_n\cap 
A_{n-1}'$ is a masa of $A_n\cap A_{n-1}'$. Also $\delta+\ad\,ih$ 
vanishes on $C$.  Let $\beta$ be the one-parameter automorphism 
group generated by $\delta+\ad\,ih$. Then there is an 
$\alpha$-cocycle $u$ such that $\beta_t=\Ad\,u_t\,\alpha_t$.  If 
$\alpha_1=\id$, then it follows that $\beta_1=\Ad\,u_1$, i.e., 
$u_1\in C$. Since $C$ is AF , we find a self-adjoint $k\in C$ such 
that $e^{ik}=u_1$. Since $(\delta+\ad\,ih)(k)=0$, one can conclude 
that $\delta+\ad(ih-ik)$ generates a one-parameter automorphism 
group $\gamma$ with $\gamma_1=\id$ such that $\gamma$ leaves each 
$A_n$ invariant.  Since $\alpha$ and $\gamma$ can be considered as 
actions of $\T$ and $\gamma$ is a cocycle-perturbation of 
$\alpha$, it follows that the crossed products 
$A\times_{\alpha}\T$ and $A\times_{\gamma}\T$ are isomorphic. 
Since $A\times_{\gamma}\T$ is AF as the inductive limit of 
$A_n\otimes C_0(\Z)$ and $A^{\alpha}$ is a hereditary 
$C^*$-subalgebra of $A\times_{\alpha}\T$, $A^{\alpha}$ must be AF.
This contradiction shows that $\D(\delta)$ cannot contain a 
canonical AF masa. \end{pf}

If $\alpha$ is not periodic in the proof of \ref{3D}, we could 
still use the following property for $\delta_{\alpha}$:
 
\medskip
\noindent{\it Condition} For any $\epsilon>0$ there exists a 
$\nu>0$ with the following property: If $u\in\D(\delta_{\alpha})$ 
is a unitary with $\|\delta_{\alpha}(u)\|<\nu$ there is a 
continuous path $(u_t)$ of unitaries in $\D(\delta_{\alpha})$ such 
that $u_0=1,\ u_1=u$, and $\|\delta_{\alpha}(u_t)\|<\epsilon,\ 
t\in [0,1]$. 

\begin{prop}\label{3F}
Let $\alpha$ be a one-parameter automorphism group of a unital 
simple AF \cstar. If $\D(\delta_{\alpha})$ contains a canonical AF 
masa, then the above Condition for $\delta_{\alpha}$ is satisfied. 
\end{prop} 
\begin{pf} 
First suppose that $A$ is finite-dimensional. Then there is an 
$h=h^*\in A$ such that $\delta_{\alpha}=\ad\,ih$. The condition 
$\|\delta_{\alpha}(u)\|<\nu$ reads
 $$ \|h-uhu^*\|<\nu.
 $$
Then the {\it Condition} follows from Theorem \ref{4A}, which will 
be given later. Note that here the choice of $\nu$ does not depend 
on $A$ nor $\delta_{\alpha}$.  

Let $(A_n)$ be an increasing sequence of finite-dimensional 
$C^*$-subalgebras of $A$ with dense union such that 
$\alpha_t(A_n)=A_n$. Let $u\in\D(\delta_{\alpha})$ be a unitary 
with $\|\delta_{\alpha}(u)\|<\nu$. Since $\cup_nA_n$ is dense in 
$\D(\delta_{\alpha})$, there is a sequence $(u_n)$ in $\cup_nA_n$ 
such that $\|u-u_n\|_{\delta_{\alpha}}\ra0$. Since 
$u_nu^*\approx1$ and $\delta_{\alpha}(u_nu^*)\approx0$, we can 
find a continuous path $(u_n(t))$ in $\D(\delta_{\alpha})$ such 
that $u_n(0)=u_n,\ u_n(1)=u$, and $\|\delta_{\alpha}(u_n(t))\|$ is 
of the order of $\|\delta_{\alpha}(u)\|$. Thus we can suppose that 
$u\in \cup_nA_n$ and the assertion follows from the previous 
paragraph. 

If $\delta_{\alpha}=\delta_{\beta}+\ad\,ih$ with 
$\beta_t(A_n)=A_n$, there is a sequence $(h_n)$ with $h_n=h_n^*\in 
A_n$ such that $\|h-h_n\|\ra0$. Then 
$\delta_{\beta}+\ad\,ih_n=\delta_{\alpha}+\ad(ih_n-ih)$ generates 
an AF locally representable action. Thus we may as well assume 
that $\alpha$ is AF locally representable. This completes the 
proof. 
\end{pf} 

If $\alpha$ is periodic and $K_1(A^{\alpha})$ is not trivial, then 
the {\em Condition} is not satisfied. But we note:
  
\begin{rem}
In the above proposition the converse does not hold. In fact the 
example in the proof of \ref{C} satisfies the above {\em 
Condition}. \end{rem}
  
By using Proposition \ref{3F} we can give more examples with the 
property that that $\D(\delta_{\alpha})$ contains no canonical AF 
masa. For example, as in the proof of Proposition \ref{3D}, 
suppose that we express $A$ as $\ol{\cup_nB_n}$ with 
$B_n=A_n\otimes C(\T)$ and that we define an $\alpha$ by defining 
$h_{n,ij}$. This time we choose $h_{n,ij}$ to be of the form:
 $$
 h_{n,ij}=\left(\begin{array}{cc}a_n&0\\0&a_n\end{array}\right)\oplus  
          \left(\begin{array}{cc}-a_n&0\\0&-a_n\end{array}\right)\oplus0\oplus\cdots\oplus0,
 $$
where $(a_n)$ is an arbitrary sequence such that $a=\inf a_n>0$.  
If we had a continuous path $(u_t)$ of unitaries in $B_n$ such 
that $u_0=1, \ u_1=z_1$, and $\|\delta_{\alpha}(u_t)\|<a$ for the 
canonical unitary $z_1$, we could reach a contradiction as 
follows. Let $H_n=\sum_j \lambda_jE_j$ be the spectral 
decomposition with $\lambda_1>\lambda_2>\cdots$. By the assumption 
we have that $\lambda_1-\lambda_2\geq a$. Since $\|[H_n,u_t]\|<a$, 
we can estimate 
 $$
 \|\sum_j (\lambda_1-\lambda_j)E_1u_tE_j\|<a,
 $$
which shows that
 $$\|E_1u_t(1-E_1)\|<1.
 $$
Since $E_1u_tE_1u_t^*E_1+\|E_1u_t(1-E_1)\|^2E_1\geq E_1$, it 
follows that $E_1u_tE_1$ is invertible. Since $E_1u_0E_1=E_1$ and 
$E_1u_1E_1=z_1E_1$ is a unitary with non-trivial $K_1$, this is a 
contradiction.  If we have a continuous path of unitaries in 
$\D(\delta_{\alpha})$ with the above property, we approximate the 
path by a path in $\cup_nB_n$ to reach the contradiction. Since 
$\cup_nB_n$ is dense in the Banach *-algebra 
$\D(\delta_{\alpha})$, this is possible. Thus we have shown that 
$\D(\delta_{\alpha})$ contains no canonical AF masa. If $\limsup 
a_n<\infty$, one can also show that $\D(\delta_{\alpha})$ is AF. 

There is a standard way to construct a one-parameter automorphism 
group $\alpha$ of a certain UHF algebra through an interaction of 
a quantum spin system \cite{BR}. If the interaction is {\em 
quantum}, we expect that any inner perturbation of $\alpha$ is not 
AF locally representable. We also expect that the quasi-free 
one-parameter automorphism group of the CAR algebra induced by a 
one-particle Hamiltonian  with continuous spectrum \cite{BR,Sak} 
or any inner perturbation of it is not AF locally representable. 
We conclude this section by posing: 

\medskip
\noindent {\em Problem} Prove the above conjecture. 

\section{A homotopy lemma}
We prove here a technical lemma which is used in the proof of 
Proposition \ref{3F}.  With an additional assumption on $h$ below 
(saying the norm is less than 1), this follows from Lemma 5.1 of 
\cite{BEEK3}. To remove this assumption we have to replace a 
certain approximation argument used there by a {\em constructive} 
argument, which will constitute the main part of the proof.

\begin{theo}\label{4A}
For any $\epsilon>0$ there exists a $\nu>0$ satisfying the 
following condition: For any unital AF algebra $A$ and $u,h\in A$ 
such that $u^*u=uu^*=1,\ h^*=h$, and $\|[h,u]\|<\nu$, there is a 
rectifiable path $(u_t)_{t\in [0,1]}$ in the unitary group of $A$ 
such that $u_0=1,\ u_1=u,\ \|[h,u_t]\|<\epsilon$, and the length 
of $(u_t)$ is smaller than $3\pi+\epsilon$. 
\end{theo}
\begin{pf}
We may assume that $A$ is finite-dimensional; in particular we 
assume that $h$ is diagonal.

Let $\delta>0$ be a sufficiently small number, which will be 
chosen later depending on $\epsilon$.

Let $f$ be a $C^{\infty}$-function on $\R$ such that $f\geq0,\ 
\int f(t)d t=1$, and $\supp\,\hat{f}\subset (-\delta,\delta)$. 
Define
 $$ x=\int f(t)e^{ith}ue^{-ith}dt.
 $$
Then it follows that $\|x\|\leq1$, $\|[h,x]\|\leq \|[h,u]\|$, and
 $$
 \|x-u\|\leq \|\int f(t)(\Ad\,e^{ith}(u)-u)dt\|\leq\int f(t)|t|dt \|[h,u]\|,
 $$
where we have used that
 $$
 \|\Ad\,e^{ith}(u)-u\|=\|\int_0^te^{ish}[ih,u]e^{-ish}ds\|\leq |t|\|[h,u]\|.
 $$
If we denote by $E_h$ the spectral measure of $h$, then we have 
that for $x^{\#}=x$ or $x^*$ and $t\in\R$,
 $$E_h(-\infty,t)x^{\#}E_h[t+\delta,\infty)=0.
 $$
We define a projection $e_k$ for each $k\in\Z$ by
 $$
 e_k=E_h[2k\delta,2(k+1)\delta).
 $$
Then there are only a finite number of non-zero $e_k$. It follows 
that $\sum_ke_k=1$ and
 $$
 xe_kx^*\leq E_h[(2k-1)\delta, (2k+3)\delta).
 $$
 
We suppose that $\|x-u\|<\mu$, where $\mu$ can be made arbitrarily 
small by choosing $\nu$ small. Since $0\leq 1-x^*x<2\mu$ and 
$0\leq xe_kx^*-(xe_kx^*)^2\leq 2\mu xe_kx^*$, we have that $0\leq 
1-\sum_k xe_kx^*\leq 2\mu$, and $\Sp(xe_kx^*)\subset \{0\}\cup 
[1-2\mu,1]$. If $x$ were a unitary (and so $xe_kx^*$ were a 
projection), we could skip most of the arguments below. What we 
will do next is to construct a unitary $v$ by using $x$ such that 
$v$ is close to $u$ and satisfies that $ve_kv^*\leq 
E_h[(2k-1)\delta,(2k+3)\delta)$. 

By multiplying $F_j=E_h[(2j-1)\delta,\infty)$ with 
$1-\sum_kxe_kx^*$ whose norm is less than $2\mu$, we get that 
 $$
 \|\sum_{k\geq j}xe_kx^*+ F_jxe_{j-1}x^*-F_j\|<2\mu,
 $$
which implies that
 $$ \|F_jxe_{j-1}x^*-xe_{j-1}x^*F_j\|<4\mu.
 $$
Since
 $$ \|(F_jx e_{j-1}x^*F_j)^2-F_jxe_{j-1}x^*F_j\|<4\mu+
  \|F_j((xe_{j-1}x^*)^2-xe_{j-1}x^*)F_j\|<6\mu,
  $$
$F_jxe_{j-1}x^*F_j$ is close to a projection for a small $\mu$. If 
we denote by $f_{j-1}^+$ the support projection of 
$F_jxe_{j-1}x^*F_j$, then we have that
 $$
 \|f_{j-1}^+-
 E_h[(2j-1)\delta,\infty)xe_{j-1}x^*E_h[(2j-1)\delta,\infty)\|<6\mu',
 $$
where $\mu'=(1-\sqrt{1-24\mu})/12\approx \mu$, which we again 
denote by $\mu$ below. Note that 
$\|f_{j-1}^+-f_{j-1}^+xe_{j-1}x^*\|<10\mu$ and that $f_{j-1}^+\leq 
E_h[(2j-1)\delta,(2j+1)\delta)$.  In the same way we denote by 
$f_{j}^-$ the support projection of 
 \BE&& E_h(-\infty,(2j+1)\delta)xe_{j}x^*E_h(-\infty,(2j+1)\delta)\\
 &&=E_h[(2j-1)\delta,(2j+1)\delta)xe_jx^*E_h[(2j-1)\delta,(2j+1)\delta);
 \EE 
then we have  that
 $$    
 \|f_{j}^--E_h[(2j-1)\delta,(2j+1)\delta)xe_{j}x^*E_h[(2j-1)\delta,(2j+1)\delta)\|<6\mu.
 $$ 
Let $f_{j}=f_{j}^-+f_{j}^+$. Then summing up the above 
calculations, we obtain that 
 \BE
 \|f_{j}-xe_{j}x^*\|&=& \|f_j-(E_h[(2j-1)\delta,(2j+1))+E_h[(2j+1)\delta,(2j+3)\delta))
  xe_jx^*\|\\
 &<& 8\mu+\|f_j-E_h[(2j-1)\delta,(2j+1)\delta)xe_jx^*E_h[(2j-1)\delta,(2j+1)\delta)\\
  &&-
  E_h[(2j+1)\delta,(2j+3)\delta)xe_jx^* E_h[(2j+1)\delta,(2j+3)\delta)\|\\
  &<& 14\mu.
 \EE
Hence if $\mu$ is small, $f_jxe_j(e_jx^*f_jxe_j)^{-1/2}$ defines a 
partial isometry with initial projection $e_j$ and final 
projection $f_j$. 

Let $g_j^-=E_h[(2j-1)\delta,(2j+1)\delta)-f_{j-1}^+$, $g_j^+=
E_h[(2j+1)\delta,(2j+3)\delta)-f_{j+1}^-$, and $g_j=g_j^-+     
g_j^+$.   Since $\|g_jxe_{j-1}x^*\|=\|g_j^-xe_{j-1}x^*\|= 
\|(1-f_{j-1}^+)E_h[(2j-1)\delta,(2j+1)\delta)xe_{j-1}x^*\|<4\mu$ 
etc.  we obtain that 
 $$
 \| g_jxe_jx^*-g_j\|<10\mu.
 $$ 
Since 
 \BE
 \|f_{j-1}xe_jx^*\|&<&6\mu+\|E_h[(2j-1)\delta,(2j+1)\delta)xe_{j-1}x^*
 E_h[(2j-1)\delta,(2j+1)\delta)xe_jx^*\|\\
 &<& 10\mu+\|E_h[(2j-1)\delta,(2j+1)\delta)xe_{j-1}x^*xe_jx^*\|\\
 &<& 12\mu,
 \EE
we have that
 \BE
 \|xe_jx^*-g_j\|&\leq&\|g_jxe_jx^*-g_j\|+\|f_{j-1}^+xe_jx^*\|
 +\|f_{j+1}^-xe_jx^*\|\\
 &<& 34\mu.
 \EE
Let
 $$
 v=\sum_jf_{2j}xe_{2j}(e_{2j}x^*f_{2j}xe_{2j})^{-1/2}+
    \sum_jg_{2j-1}xe_{2j-1}(e_{2j-1}x^*g_{2j-1}xe_{2j-1})^{-1/2},
 $$                                                                
which is the unitary part of the polar decomposition of 
$y=\sum_jf_{2j}xe_{2j}+\sum_jg_{2j-1}xe_{2j-1}$. Since $0\leq 
1-yy^*<14\mu$, we have that $\|v-y\|<1/\sqrt{1-14\mu}-1$. Note 
also that $ ve_jv^*\leq E_h[(2j-1)\delta,(2j+3)\delta)$. Since 
$\|\sum(f_{2j}-1)xe_{2j}\|^2=\sup_j 
\|(f_{2j}-1)xe_{2j}x^*(f_{2j}-1)\|\leq 
\sup_j\|xe_{2j}x^*-f_{2j}\|$,  we get that 
 $$
 \|\sum(f_{2j}-1)xe_{2j}\|<\sqrt{14\mu}.
 $$
Since $\|(g_{2j-1}-1)xe_{2j-1}x^*(g_{2j-1}-1)\|<34\mu$, we get 
that 
 $$
 \|\sum(g_{2j-1}-1)xe_{2j-1}\|<\sqrt{34\mu}.
 $$
 Since 
$\|y-x\|\leq 
\|\sum(f_{2j}-1)xe_{2j}\|+\|\sum(g_{2j-1}-1)xe_{2j-1}\|$, we get 
$\|y-x\|<\sqrt{14\mu}+\sqrt{34\mu}<10\sqrt{\mu}$. Hence we get 
that if $\mu$ is sufficiently small,
 $$
 \|v-u\|<\|v-y\| +\|y-x\|+\|x-u\|<1/\sqrt{1-14\mu}-1+10\sqrt{\mu}+\mu
 <10(\mu+\sqrt{\mu}).
 $$
We assume that the constant $10(\mu+\sqrt{\mu})$ is sufficiently 
small. 

Let $k=\sum_j 2j\delta E_h[2j\delta,2(j+1)\delta)=\sum 2j\delta 
e_j$. Then $\|h-k\|<2\delta$ and $\|[k,u]\|\leq 2\|h-k\|+\|[h,u]\|
<4\delta+\|[h,u]\|<4\delta+\nu$. Since $ve_kv^*\leq 
e_{k-1}+e_k+e_{k+1}$, we have that $k-2\delta\leq vkv^*\leq 
k+2\delta$. Hence it follows that  $\|vkv^*-k\|\leq 2\delta$. 
Since $\|[vu^*,k]\|\leq \|[v,k]\|+\|[u,k]\|<6\delta+\nu$ and 
$vu^*=e^{ia}$ with $a^*=a\approx 0$, we get that 
$\|[a,k]\|\approx0$ (up to the order of $6\delta+\nu$). We take a 
continuous path $t\in[0,1]\mapsto w_t=e^{ita}u$ of length $\|a\|$. 
Then since $[k,w_t]=[k,e^{ita}]u+e^{ita}[k,u]\approx0$ (up to the 
order of $10\delta+2\nu$) and $w_1=v$, we may replace $u$ by $v$. 
From now on we can proceed as in the proof of Lemma 5.1 of 
\cite{BEEK3}. 

Let $E_n=\sum_{j\geq n}e_j$. Then $k$ equals $2\delta \sum_{n>m} 
E_n+2m\delta$, where $m$ is the biggest integer satisfying 
$E_{m}=1$, and the sequence $(E_n)_{n\geq m}$ of projections 
decreases from 1 to 0 as $n$ increases. Let $F_n=vE_nv^*$. Then we 
have that $E_{n+1}\leq F_n\leq E_{n-1}$. Since $F_{2n+2}\leq 
F_{2n+1}\leq F_{2n}$ and $F_{2n+2}\leq E_{2n+1}\leq F_{2n}$, we 
find a continuous path $(w_t)$ of unitaries of length at most 
$\pi$ such that $w_0=1$, $[w_t,F_{2n}-F_{2n+2}]=0$, and 
$w_1(F_{2n+1}-F_{2n+2})w_1^*=E_{2n+1}-F_{2n+2}$ for all $n$. Since 
$\|vkv^*-4\delta\sum_{2n>m} F_{2n}-2m\delta\|\leq 2\delta$, we 
have that $\|w_tvkv^*w_t^*-vkv^*\|\leq 4\delta$ and hence 
$\|[w_tv,k]\|\leq 6\delta$. Next we find a continuous path $(z_t)$ 
of unitaries of length at most $\pi$ such that $z_0=1$, 
$[z_t,E_{2n-1}-E_{2n+1}]=0$, and 
$z_1(F_{2n}-E_{2n+1})z_1^*=E_{2n}-E_{2n+1}$. Since 
$w_1vkv^*w_1^*=2\delta(\sum_{2n>m}F_{2n}+\sum_{2n+1>m}E_{2n+1})+2m\delta$ 
and  $\|w_1vkv^*w_1^*-4\delta\sum_{2n+1>m} E_{2n+1}-2m\delta\|\leq 
2\delta$, we get that $\|z_tw_1vkv^*w_1^*z_t-w_1vkv^*w_1^*\|\leq 
4\delta$ and hence $\|[z_tw_1v,k]\|\leq 10\delta$. Since 
$z_1w_1vkv^*w_1^*z_1^*=k$, we can find a continuous path of 
unitaries from $z_1w_1v$ to 1 in the commutant of $k$, whose 
length is at most $\pi$. (Here we use the fact that the unitary 
group of $eAe$ for any projection $e\in A$ is connected.) Note 
that the path obtained by combining these three paths has length 
at most $3\pi$. 

The above calculations show that we can choose $\delta$ just 
depending on $\epsilon$. (For example $\delta$ should be smaller 
than $\epsilon/15$ and much smaller than 1.) Then we choose $\nu$ 
independently (such that  $\nu$ is smaller than $\epsilon/30$ and
$10(\mu+\sqrt{\mu})$ is much smaller than $\epsilon$, where $\mu$ 
is proportional to $\nu$ as shown at the beginning of the proof).   
This concludes the proof.

\end{pf}

\medskip
\small

\end{document}